\documentclass[12pt]{article}
\usepackage{amsfonts}

\usepackage{latexsym}
\usepackage{amssymb}
\usepackage{epsfig}
\usepackage{amsmath}
\usepackage{amsthm}
\usepackage[mathscr]{eucal}
\usepackage{subfigure}
\usepackage{graphicx}
\usepackage{hyperref}

\newtheorem{Lemma}{Lemma}

\newtheorem{Theorem}[Lemma]{Theorem}

\newtheorem{Corollary}[Lemma]{Corollary}

\renewcommand{\qed}{\hfill{\ \ \rule{2mm}{2mm}} \vspace{0.2in}}

\newcommand{\ind}{1\hspace{-2.3mm}{1}}

\begin{document}

\title{Records from partial comparisons and discrete approximations}
\author{ \textbf{Ghurumuruhan Ganesan}
\thanks{E-Mail: \texttt{gganesan82@gmail.com} } \\
\ \\
New York University, Abu Dhabi }
\date{}
\maketitle

\begin{abstract}
In this paper we study records obtained from partial comparisons within a sequence of independent
and identically distributed (i.i.d.) random variables, indexed by positive integers, with a common density~\(f.\)
Our main result is that if the comparison sets along a subsequence of the indices
satisfy a certain compatibility property, then the corresponding record events are independent.
Moreover, the record event probabilities do not depend on the density~\(f\) and we obtain closed form
expressions for the distribution of~\(r^{th}\) record value for any integer~\(r \geq 1.\) 

Our proof techniques extend to the discrete case as well and we estimate the difference in record event probabilities
associated with a continuous random variable~\(X\)
and its discrete approximations.



\vspace{0.1in} \noindent \textbf{Key words:} Records, partial comparisons, discrete approximations.

\vspace{0.1in} \noindent \textbf{AMS 2000 Subject Classification:} Primary:
60J10, 60K35; Secondary: 60C05, 62E10, 90B15, 91D30.
\end{abstract}

\bigskip

\setcounter{equation}{0}
\renewcommand\theequation{\thesection.\arabic{equation}}
\section{Introduction}\label{intro}
Records in a sequence of independent and identically distributed (i.i.d.) random variables with a common density~\(f\) have been extensively studied in the context of statistical estimation (see Ahsanullah and Nevzorov (2015), Balakrishnan and Chan (1998), Nevzorov (1986) and references therein). When the random variables are continuous valued, the probability that the~\(k^{th}\) random variable is a record grows inversely in the index~\(k\) and does not depend on the distribution.
Moreover, the record events are mutually independent.

The main feature of the above setup is that the records are determined by comparing the value of the random variable at a certain index~\(k\) with the values of \emph{all} previous indices. In this paper, we study records determined from partial comparisons where we determine whether or not a record occurs at index~\(k\) by considering previous values belonging to a \emph{subset} of the indices~\(\{1,2,\ldots,k-1\}.\) Moreover, we consider records with indices belonging to an arbitrary subsequence of the original sequence. Our main result (see Theorem~\ref{prop1}) is that if the subsequence and the corresponding collection of comparison sets are compatible in a certain sense, then the record events along the subsequence are independent and the corresponding probabilities do not depend on the density~\(f.\) We also obtain closed form expressions for the value of the~\(r^{th}\) record for any integer~\(r \geq 1.\)


Our proof techniques extend to the discrete case as well and allows us to estimate record properties associated with a continuous random variable~\(X\) using a sequence of discrete approximations that converge weakly to~\(X\) (see discussion prior to Theorem~\ref{t1t2_thm} for more details).

\subsection*{Records from partial comparisons}
Let~\(\{X_i\}_{1 \leq i \leq n}\) be a sequence of independent and identically distributed (i.i.d) random variables with density~\(f(.)\) and cumulative distribution function (cdf)
\begin{equation}\label{cdf_def}
F(y) := \mathbb{P}(X_1 < y) = \int_{0}^{y} f(z) dz.
\end{equation}
The random variables~\(\{X_i\}_{1 \leq i \leq n}\) are defined on the probability space~\((\Omega,{\cal F},\mathbb{P})\)
and we assume throughout that all densities are continuous. 
For every integer~\(n \geq 1\) we associate a comparison set~\({\cal C}(n) \subseteq \{1,2,\ldots n-1\}\) and say that a record occurs at index~\(n\) if the event
\begin{equation}\label{ak_def}
A_n := \left\{X_n > \max_{j \in {\cal C}(n)} X_j\right\}
\end{equation}
occurs. We say that~\(A_n\) is the record event occurring from partial comparison using the set~\({\cal C}(n).\) If~\({\cal C}(n) = \{1,2,\ldots,n-1\},\) then~\(A_n\) is the usual record event occurring from total comparison with all previous values.

Apart from partial comparisons, we are also interested in studying records from a subsequence~\(\{n_t\}_{ t \geq 1}\) of the time index~\(\{n\}.\) Let~\({\cal I} := \{n_t\}\) be any strictly increasing sequence of integers and let~\({\cal C} := \{{\cal C}(n_t)\}\) be the corresponding collection of comparison sets. We assume throughout~\({\cal I}\) and~\({\cal C}\) are compatible in the sense that the following two conditions hold:\\
\((a1)\) The comparison sets~\({\cal C}(n_{1}) \subset {\cal C}(n_2) \ldots \) are strictly increasing and\\
\((a2)\) For all~\(t \geq 2,\) the index~\(n_{t-1}  \in {\cal C}(n_{t}).\)\\
For example, let~\(n_1 = 1, \{n_t\}_{t\geq 2}\) be any strictly increasing sequence and let~\({\cal C}(1) = \{1\}\) and for~\(t \geq 2\) let \[{\cal C}(n_t) := \{1\} \bigcup \bigcup_{u=1}^{t-1} \{n_u\}.\] Conditions~\((a1)-(a2)\) are satisfied in this case.


Compatibility of the subsequence~\(\{n_t\}\) and the corresponding comparison sets allow us to study record events with partial comparisons and we have the following result.
\begin{Theorem}~\label{prop1} For any integer~\(r \geq 1\) the probability of a record occurring at index~\(n_r\) is
\begin{equation}\label{vnr}
\mathbb{P}(A_{n_r}) = \frac{1}{c(n_r)}  := \frac{1}{\#{\cal C}(n_r)+1}.
\end{equation}
For any finite set of indices~\(\{i_1,\ldots,i_k\} \subset \{n_t\},\) the corresponding record events~\(\{A_{i_j}\}_{1 \leq j \leq k}\) are mutually independent in the sense that
\begin{equation}\label{mut_ind}
\mathbb{P}\left(A_{i_1} \bigcap A_{i_2} \bigcap \ldots \bigcap A_{i_k}\right) = \prod_{j=1}^{k} \mathbb{P}\left(A_{i_j}\right) = \frac{1}{c(i_1)\cdot c(i_2) \cdots c(i_k)}.
\end{equation}
\end{Theorem}
In Section~\ref{pfs1}, we use induction to prove~(\ref{mut_ind}).

Using Theorem~\ref{prop1}, we now study records along the subsequence~\(\{n_t\}\) by letting~\(j \geq 1\) and defining
\begin{equation}\label{rec}
R_j := \sum_{1 \leq t \leq j} \ind(A_{n_t})
\end{equation}
to be the number of records in the first~\(j\) indices of the subsequence. We have the following result.
\begin{Theorem}\label{cor_samp}  For any~\(j \geq 2\)
\begin{equation}\label{pk_est1_samp}
\mathbb{E}R_j = I_j := \sum_{1 \leq t \leq j} \frac{1}{c(n_t)} \text{ and } var(R_j) = I_j - \sum_{1 \leq t \leq j}\frac{1}{c^2(n_t)}.
\end{equation}
If~\(I_j \longrightarrow \infty\) as~\(j \rightarrow \infty,\)
then~\(\frac{R_j}{I_j} \longrightarrow 1\) a.s.\ as~\(j \rightarrow \infty.\)
\end{Theorem}
Thus the number of records in the first~\(j\) indices of the sequence~\(\{n_t\}\) roughly grows as~\(I_j.\) To compute the individual record values, define~\(L(0) = 0\) and for~\(r \geq 1\) the~\(r^{th}\) record time
\begin{equation}\label{lr_def}
L(r) := \inf \left\{t \geq 1: n_t \geq L(r-1)+1 \text{ and } A_{n_t} \text{ occurs }\right\}
\end{equation}
be the index where the~\(r^{th}\) record occurs along the sequence~\(\{n_t\}.\)
\begin{Theorem}\label{thm_rec_val} For~\(r \geq 1\) and~\(x > 0\) we have that
\begin{equation}\label{xlr_x}
\mathbb{P}\left(X_{L(r)} < x\right) = \sum_{t \geq r} F^{n_t}(x) \mathbb{P}\left(L(r) = n_t\right) = \mathbb{E}F^{L(r)}(x).
\end{equation}
\end{Theorem}

So far, we have studied records obtained from partial comparisons within a subsequence of the original sequence of random variables. If on the other hand we set~\(n_t = t\) and~\({\cal C}(n_t) = \{1,2,\ldots,t-1\}\) for all~\(t,\) then we recover the usual definition of records obtained from total comparisons. In this case~\(c(t) = \#{\cal C}(t)+1 = t\) for all~\(t\) and so~\(I_j = H_j := \sum_{t=1}^{j} \frac{1}{t}.\) Moreover, there is a constant~\(\gamma > 0\) such that~\(H_j -\log{j} \longrightarrow \gamma\) as~\(j \rightarrow \infty\) and so~\(\frac{I_j}{\log{j}} \longrightarrow 1\) as~\(j \rightarrow \infty.\) Consequently, from Theorem~\ref{cor1} we get the following result.
\begin{Corollary}\label{cor1}
For any~\(j \geq 2\)
\begin{equation}\label{pk_est1}
\mathbb{P}(A_j) = \frac{1}{j}, \mathbb{E}R_j = H_j,  var(R_j) = H_j - \sum_{t=1}^{j}\frac{1}{t^2}
\end{equation}
and~\(\frac{R_j}{\log{j}} \longrightarrow 1\) a.s.\ as~\(j \rightarrow \infty.\)
\end{Corollary}

\subsection*{Records from discrete approximations}
Suppose now that~\(X\) is a continuous random variable with a continuously differentiable density~\(f(.)\) and suppose there exists a~\(M > 0\) such that\\\(f(x) = 0\) for all~\(x < 0\) and~\(x >M.\)  For~\(m \geq 1,\) let~\(Y_m\) be a discrete valued random variable with distribution
\begin{equation}\label{yn_dist}
\mathbb{P}_m\left(Y_m = l\right) = \frac{f\left(\frac{l}{m}\right)}{\sum_{l = 0}^{Mm} f\left(\frac{l}{m}\right)} ,\;\;l=0,1,2,\ldots, Mm
\end{equation}
and let~\(X_{m,1},\ldots,X_{m,n}\) be independent and identically distributed (i.i.d)
with same distribution as~\(Y_m,\) defined on the probability space~\((\Omega_{m}, {\cal F}_m, \mathbb{P}_m).\)
As in~(\ref{ak_def}), we let~\(A_{n} =  \left\{X_{m,n} > \max_{j \in {\cal C}(n)} X_{m,j}\right\}\) denote
the event that a record occurs at index~\(n.\)

By definition, the sequence of random variables~\(\{Y_m\}\) converge weakly to~\(X\)
and so for any fixed index~\(i,\) we use~(\ref{vnr}) to obtain that the corresponding record event probability~\[\mathbb{P}_m(A_{i}) \longrightarrow \frac{1}{c(i)} = \frac{1}{\#{\cal C}(i)+1},\] as~\(m \rightarrow \infty.\) The result below provides estimates on difference between the probabilities of the record events  and their corresponding asymptotic values.
\begin{Theorem}\label{t1t2_thm} Suppose there is a constant~\(C > 0\) such that
\begin{equation}\label{f_bd}
\max(|f(x)|, |f'(x)|) \leq C \text{ for all } 0 \leq x \leq M.
\end{equation}
For any fixed finite set of indices~\(\{i_1,\ldots,i_k\} \subset \{n_t\}\) there is a constant~\(D = D(i_1,\ldots,i_k) > 0\) such that
\begin{equation}\label{mut_ind_disc}
\left|\mathbb{P}_m\left(A_{i_1} \bigcap A_{i_2} \bigcap \ldots \bigcap A_{i_k}\right) - \frac{1}{c(i_1)\cdot c(i_2) \cdots c(i_k)}\right|
\leq \frac{D}{m}
\end{equation}
for all~\(m\) large.
\end{Theorem}
Thus the record events~\(\{A_{i_j}\}_{1 \leq j \leq k}\) are asymptotically independent and converge to the distribution-free value determined in~(\ref{vnr}), as~\(m \rightarrow \infty.\) If we think of the random variable~\(Y_m\) as discrete approximations of the random variable~\(X,\) then Theorem~\ref{t1t2_thm} estimates the difference between record event probabilities associated with~\(X\) and its discrete approximations, in terms of the approximation interval~\(\frac{1}{m}.\) Relation~(\ref{mut_ind_disc}) can then be used to estimate any record properties associated with~\(X,\)
using the discrete approximations~\(\{Y_m\}.\)

The paper is organized as follows. In Section~\ref{pfs1}, we prove Theorems~\ref{prop1},~\ref{cor_samp} and~\ref{thm_rec_val} and in Section~\ref{pfs2}, we prove Theorem~\ref{t1t2_thm}.


\setcounter{equation}{0}
\renewcommand\theequation{\arabic{section}.\arabic{equation}}
\section{Proof of Theorems~\ref{prop1},~\ref{cor_samp} and~\ref{thm_rec_val}}\label{pfs1}
We first prove~(\ref{mut_ind}) for~\(k=1\) and~\(k=2\) and obtain the general result by induction.
Throughout we use the following relations: Suppose~\(Y\) and~\(Z\) are independent random variables with continuous cdfs and suppose~\(A\) is an event independent of~\(Z.\) If~\(Z\) has a continuous density~\(f(.),\) then applying Fubini's theorem we have that
\begin{equation}\label{x1x2}
\mathbb{P}\left(A \bigcap \{Y < Z\}\right) = \int_{0}^{\infty} \mathbb{P}\left(A \bigcap \{Y < z\}\right) f(z) dz,
\end{equation}
For notational simplicity, we use the phrase ``conditioning on~\(Z=z\)" whenever we refer to~(\ref{x1x2}) or an analogous estimate
obtained via Fubini's theorem. Also if~\(g\) is a function with a continuous derivative~\(g',\) then by the fundamental theorem of calculus we have for~\(0 \leq a \leq b \leq \infty\) that
\begin{equation}\label{fund_calc}
\int_{a}^{b} g'(x) dx = g(b) - g(a).
\end{equation}
We refer to~(\ref{fund_calc}) as the definite integral property.

We now prove a slightly stronger result than~(\ref{mut_ind}) for future use. For~\(x > 0\) and integer~\(r \geq 1\) define
\begin{equation}\label{akx}
A_r(x) := A_r \bigcap \left\{X_r < x\right\} = \left\{ \max_{j \in {\cal C}(r)} X_j < X_r < x\right\}
\end{equation}
and recall that~\(c(r) = \#{\cal C}(r) + 1.\) We have that
\begin{equation}\label{k_one_case}
\mathbb{P}(A_r(x)) = \frac{F^{c(r)}(x)}{c(r)}.
\end{equation}
and the events~\(\{A_{i_j}\}_{1 \leq j \leq k-1} \cup A_{i_k}(x)\) are mutually independent in the sense that if
\begin{equation}\label{bj_def}
B_k := \bigcap_{j=1}^{k} A_{i_j} \text{ and } B_k(x) := \bigcap_{j=1}^{k-1} A_{i_j} \bigcap A_{i_k}(x),
\end{equation}
then
\begin{eqnarray}
\mathbb{P}\left(B_k(x)\right) &=& \left(\prod_{j=1}^{k-1} \mathbb{P}\left(A_{i_j}\right)\right)\mathbb{P}\left(A_{i_k}(x)\right) \nonumber\\
&=& \frac{1}{c(i_1)\cdot c(i_2) \cdots c(i_{k-1})} \frac{F^{c(i_k)}(x)}{c(i_k)} \nonumber\\
&=& \left(\prod_{j=1}^{k}\mathbb{P}\left(A_{i_j}\right)\right)F^{c(i_k)}(x) \nonumber\\
&=& \mathbb{P}(B_k)F^{c(i_k)}(x). \label{mut_ind_2}
\end{eqnarray}

\emph{Proof of~(\ref{mut_ind_2})}: We first prove~(\ref{k_one_case}). Conditioning on~\(X_{n_r} = y\) we get that
\begin{equation}\label{pr_ak_ev}
\mathbb{P}\left(A_{n_r}(x)\right) = \int_{0}^{x} \mathbb{P}\left(\max_{j \in {\cal C}(n_r)} X_j < y\right) f(y) dy = \int_{0}^{x} F^{c(n_r)-1}(y) f(y)dy
\end{equation}
and using the definite integral property~(\ref{fund_calc}) with~\(g = \frac{F^{c(n_r)}}{c(n_r)},\) the right hand side of~(\ref{pr_ak_ev}) evaluates to~\(\frac{F^{c(n_r)}(x) - F^{c(n_r)}(0)}{c(n_r)}= \frac{F^{c(n_r)}(x)}{c(n_r)},\) proving also~(\ref{mut_ind_2}) for~\(k=1.\)

To prove the induction step, we now assume that the relation~(\ref{mut_ind_2}) is true some integer~\((k-1) \geq 2\) and without loss of generality,
also assume that~\(i_1 < i_2 < \ldots <i_k.\) If the event~\(B_{k-1} := \bigcap_{j=1}^{k-1} A_{i_j}\) occurs, then since~\({\cal C}(i_{k-1}) \subset {\cal C}(i_k)\) (condition~\((a1)\)), the term~\(X_{i_k}\) is a record if and only if
\begin{equation}\label{xik1}
X_{i_k} > \max_{j \in {\cal C}(i_k) \setminus {\cal C}(i_{k-1})} X_j
\end{equation}
if and only if
\begin{equation}\label{xik}
X_{i_k} > \max_{j \in {\cal C}(i_k) \setminus \left({\cal C}(i_{k-1}) \cup i_{k-1}\right)} X_j \text{ and }X_{i_k} > X_{i_{k-1}},
\end{equation}
since~\(i_{k-1} \in {\cal C}(i_k)\) (condition~\((a2)\)). Therefore~\(\ind\left(B_k(x)\right) = \ind\left(A_{i_k}(x) \bigcap B_{k-1}\right)\) equals
\begin{equation}\label{ai2_ai1_k}
\ind\left(x > X_{i_k} > \max_{j \in {\cal C}(i_k) \setminus \left({\cal C}(i_{k-1}) \cup i_{k-1}\right)} X_j\right) \ind(x > X_{i_k} > X_{i_{k-1}}) \ind(B_{k-1}).
\end{equation}
The event~\( \ind(B_{k-1})\) depends only on the values of~\(\{X_j\}_{1 \leq j \leq i_{k-1}}\) and so conditioning on~\(X_{i_k}=z\)
we get from~(\ref{ai2_ai1_k}) that
\begin{equation}\label{pa2a1}
\mathbb{P}\left(B_k(x)\right) = \int_{z=0}^{x} t_{i_1,i_2,\ldots,i_k}(z)f(z) dz,
\end{equation}
where~
\begin{eqnarray}
t_{i_1,i_2,\ldots,i_k}(z) &=&  \mathbb{P}\left(\max_{j \in {\cal C}(i_k) \setminus \left({\cal C}(i_{k-1}) \cup i_{k-1}\right)} X_j < z\right) \mathbb{P}\left(B_{k-1} \bigcap \left\{ X_{i_{k-1}} < z\right\} \right) \nonumber\\
&=& F^{c(i_k)-c(i_{k-1})-1}(z)\mathbb{P}\left(B_{k-1} \bigcap \left\{ X_{i_{k-1}} < z\right\} \right) \nonumber\\
&=& F^{c(i_k)-c(i_{k-1})-1}(z)\mathbb{P}\left(B_{k-1}(z)\right) \nonumber
\end{eqnarray}
so that
\begin{equation}\label{pbb0}
\mathbb{P}\left(B_k(x)\right) = \int_{z=0}^{x} F^{c(i_k)-c(i_{k-1})-1}(z)\mathbb{P}\left(B_{k-1}(z)\right) f(z)dz.
\end{equation}
We recall from the strictly increasing comparison set sequence condition~\((a1)\) that~\(c(i_k) \geq c(i_{k-1})+1.\)

By induction assumption, we get from~(\ref{mut_ind_2}) that
\begin{equation}
\mathbb{P}\left(B_{k-1}(z)\right)  = \left(\prod_{j=1}^{k-2} \mathbb{P}\left(A_{i_j}\right)\right)\mathbb{P}\left(A_{i_{k-1}}(z)\right)
= \left(\prod_{j=1}^{k-1} \mathbb{P}\left(A_{i_j}\right)\right)F^{c(i_{k-1})}(z). \label{ind_asm}
\end{equation}
Substituting~(\ref{ind_asm}) into~(\ref{pbb0}) gives that
\begin{eqnarray}
\mathbb{P}\left(B_k(x)\right) &=&  \left(\prod_{j=1}^{k-1} \mathbb{P}\left(A_{i_j}\right) \right)\int_{z=0}^{x} F^{c(i_k)-1}(z)f(z)dz \nonumber\\
&=& \left(\prod_{j=1}^{k-1} \mathbb{P}\left(A_{i_j}\right) \right)\frac{F^{c(i_k)}(x)}{c(i_k)} \nonumber\\
&=& \left(\prod_{j=1}^{k-1} \mathbb{P}\left(A_{i_j}\right) \right)\mathbb{P}(A_{i_k}(x)) \nonumber
\end{eqnarray}
by~(\ref{k_one_case}), proving the induction step.~\(\qed\)

\emph{Proof of Theorem~\ref{cor_samp}}: The proof of the mean and variance estimates in~(\ref{pk_est1}) follow from the mutual independence relation~(\ref{mut_ind}) and the proof of the a.s.\ convergence follows by using~(\ref{mut_ind}) with~\(k=2\) together with Theorem~2.3.8 of Durrett (2013).~\(\qed\)

\emph{Proof of Theorem~\ref{thm_rec_val}}: We recall from~(\ref{akx}) that for~\(x > 0,\) the event~\(A_k(x) = A_k \bigcap \{X_k < x\}.\) Also denoting~\(m \ll n_{j}\) to be the set of all indices in~\(m \in \{n_t\}\) less or equal to~\(n_j\) we get that
\begin{equation}\label{l1_def}
\{L(1) = n_j\} = \bigcap_{u \ll n_{j-1}} A_u^c \bigcap A_{n_j}.
\end{equation}
Therefore
\begin{eqnarray}
\mathbb{P}\left(X_{L(1)} < x\right) &=& \sum_{j \geq 1} \mathbb{P}\left(\{L(1) = n_j\} \bigcap  \{X_{n_j} < x\}\right) \nonumber\\
&=& \sum_{j \geq 1} \mathbb{P}\left(\bigcap_{u  \ll n_{j-1}} A_u^c  \bigcap A_{n_j} \bigcap \{X_{n_j} < x\}\right) \nonumber\\
&=&  \sum_{j \geq 1} \mathbb{P}\left(\bigcap_{u  \ll  n_{j-1}} A_u^c \bigcap A_{n_j}(x) \right) \nonumber\\
&=& \sum_{j \geq 1} \mathbb{P}\left(\bigcap_{u  \ll n_{j-1}} A_u^c \bigcap A_{n_j}\right) F^{n_j}(x), \label{xl1}
\end{eqnarray}
using the mutual independence relation~(\ref{mut_ind_2}). From~(\ref{l1_def}) and~(\ref{xl1}) we get that
\[\mathbb{P}\left(X_{L(1)} < x\right) = \sum_{j \geq 1} \mathbb{P}\left(L(1) = n_j\right) F^{n_j}(x),\] proving~(\ref{xlr_x}) with~\(r = 1.\) The proof for general~\(r\) is analogous since
\begin{eqnarray}
\mathbb{P}(X_{L(r)}< x) &=& \sum_{j \geq r} \mathbb{P}\left(\{L(r) = n_j\} \bigcap \{X_{n_j}<x\}\right) \nonumber\\
&=&  \sum_{j \geq r} \sum \mathbb{P}\left(W(j_1,\ldots,j_{r-1},j) \bigcap \{X_{n_j}<x\}\right),  \label{xlr_c}
\end{eqnarray}
where the second summation is taken over all~\(1 \leq j_1< \ldots < j_{r-1} < j\) and the event
\[W(j_1,\ldots,j_{r-1},j) := \{L(1)=n_{j_1},L(2) = n_{j_2},\ldots,L(r-1)=n_{j_{r-1}}, L(r) = n_j\}.\]
The event~\(W(.)\) can be written as the intersection of the events~\(\{A_u\}\) as~\[W(j_1,\ldots,j_{r-1},j) = C_{n_j} \bigcap A_{n_j}\]
where~\(C_{n_j}\) equals
\begin{equation} \nonumber
\bigcap_{m_1 \ll n_{j_1-1}} A_{m_1}^c \bigcap  A_{n_{j_1}} \bigcap_{n_{j_1+1} \ll m_2 \ll n_{j_2-1}} A_{m_2}^c \bigcap A_{n_{j_2}} \bigcap \ldots \bigcap_{n_{j_{r-1}+1} \ll m_r \ll n_{j-1}} A_{m_{r}}^c.
\end{equation}
Therefore
\[W(j_1,\ldots,j_{r-1},j)  \bigcap \{X_{n_j} < x\}  = C_{n_j} \bigcap A_{n_j}(x)\]
and using the mutual independence relation~(\ref{mut_ind_2}), we therefore get that
\[\mathbb{P}\left(W(j_1,\ldots,j_{r-1},j)\bigcap \{X_{n_j}<x\} \right) = \mathbb{P}\left(W(j_1,\ldots,j_{r-1},j)\right)F^{n_j}(x) \]
and substituting back in~(\ref{xlr_c}), we get
\begin{eqnarray}
\mathbb{P}(X_{L(r)}< x) &=& \sum_{j \geq r} \sum_{1 \leq j_1 < \ldots < j_{r-1} < j} \mathbb{P}\left(W(j_1,\ldots,j_{r-1},j)\right) F^{n_j}(x) \nonumber\\
&=& \sum_{j \geq r} \mathbb{P}(L(r) = n_j)F^{n_j}(x) , \nonumber
\end{eqnarray}
proving~(\ref{xlr_x}).~\(\qed\)


\setcounter{equation}{0}
\renewcommand\theequation{\arabic{section}.\arabic{equation}}
\section{Proof of Theorem~\ref{t1t2_thm}}\label{pfs2}
For sequences~\(\{p_m(l)\}_{m \geq 1, 0 \leq l \leq Mm}\) and~\(\{q_m(l)\}_{m \geq 1, 0 \leq l \leq Mm},\) we denote
\[p_m(l) \approx q_m(l) \text{ if } \max_{0 \leq l \leq Mm} |p_m(l) - q_m(l)| \leq \frac{D}{m}\]
for all~\(m\) large and some constant~\(D > 0.\) Throughout, all sequences are positive and all constants are independent of~\(m.\)

The following Lemma is used in the proof of Theorem~\ref{t1t2_thm}.
For a constant~\(r \geq 1\) and~\(0 \leq l_1 \leq l \leq Mm,\) let
\begin{equation}\label{thet_def}
\Theta_m(l) = \Theta_m(l,r):= \frac{1}{m}\sum_{l_1=0}^{l-1} F^{r-1}\left(\frac{l_1}{m}\right)f\left(\frac{l_1}{m}\right),
\end{equation}
\begin{equation}\label{gm_def}
g_m(l) = \mathbb{P}(Y_m =l)  := \frac{f\left(\frac{l}{m}\right)}{\sum_{l = 0}^{Mm} f\left(\frac{l}{m}\right)}
\end{equation}
and let~\(G_m(l) = \sum_{l_1=0}^{l-1} g_m(l_1),\)
\begin{Lemma}\label{thet_lem}
We have that
\begin{equation}\label{gapp2}
\Theta_m(l) \approx \frac{F^{r}\left(\frac{l-1}{m}\right)}{r} \approx \frac{F^{r}\left(\frac{l}{m}\right)}{r}.
\end{equation}
Consequently,
\begin{equation}\label{gapp4}
\frac{1}{m}\sum_{l = 0}^{Mm} f\left(\frac{l}{m}\right) \approx 1, \;\;\;\;\;\;G^{r}_m(l) \approx F^{r}\left(\frac{l}{m}\right)
\end{equation}
and
\begin{equation}\label{gapp3}
\sum_{l_1=0}^{l-1} F^{r-1}\left(\frac{l_1}{m}\right) g_m(l_1) \approx \frac{1}{r}\cdot F^{r}\left(\frac{l}{m}\right)
\end{equation}
\end{Lemma}

\emph{Proof of Lemma~\ref{thet_lem}}: To prove the first relation~(\ref{gapp2}), we let~\(w(z) := F^{r-1}(z) f(z)\) where we recall that~\(F(y) = \int_{0}^{y} f(z) dz.\)  Since~\(f\) and~\(f'\) are both bounded and~\(F(z) \leq 1,\) the derivative
\begin{equation}\label{der_w}
|w'(z)| = F^{r-2}(z)|(r-1)f^2(z) + F(z) f'(z)| \leq (r-1)|f^2(z)| + |f'(z)| \leq D,
\end{equation}
for some constant~\(D  = D(r) > 0.\)
Therefore
\begin{eqnarray}
\left|\Theta_m(l) - \int_{0}^{\frac{l-1}{m}} w(z) dz \right| &=& \left| \sum_{l_1=1}^{l-1} \int_{\frac{l_1-1}{m}}^{\frac{l_1}{m}}\left(w(z) - w\left(\frac{l_1}{m}\right) \right)dz\right| \nonumber\\
&\leq& \sum_{l_1=1}^{l-1} \int_{\frac{l_1-1}{m}}^{\frac{l_1}{m}}\left|w(z) - w\left(\frac{l_1}{m}\right)\right|dz. \label{gapp}
\end{eqnarray}
Using the mean value theorem, we have some~\(z < z_{l_1} < \frac{l_1}{m}\) that
\[\left|w\left(\frac{l_1}{m}\right)-w(z)\right| = |w'(z_{l_1})| \left|\frac{l_1}{m}-z\right| \leq \frac{D}{m},\]
using~(\ref{der_w}). Therefore
\begin{equation}\label{gml1}
\left|\Theta_m(l) - \int_{0}^{\frac{l-1}{m}} w(z) dz \right| \leq \frac{D}{m}.
\end{equation}
From~(\ref{gml1}) and the fact that~\(\int_{0}^{y} w(z) dz = \frac{F^{r}(y)}{r},\) we obtain the first relation in~(\ref{gapp2}).

To prove the second relation in~(\ref{gapp2}), we use the fact that~\(|w| \leq |f|\) is bounded to get that
\begin{equation}\label{fml1}
\frac{1}{r}\left|F^{r}\left(\frac{l}{m}\right) - F^{r}\left(\frac{l-1}{m}\right) \right| = \left|\int_{\frac{l-1}{m}}^{\frac{l}{m}}w(z) dz \right| \leq \frac{C}{m}.
\end{equation}
To prove the first relation in~(\ref{gapp4}) we use~(\ref{gapp2}) with~\(r=1\) to get that
\[\frac{1}{m}\sum_{l=0}^{Mm} f\left(\frac{l}{m}\right) \approx F(M) = 1.\]

It suffices to the second relation in~(\ref{gapp4}) for~\(r=1.\) We use the fact that if~\(0 \leq a_m(l) \approx b_m(l) \leq 1\) and~\(c_m(l) \approx 1\) then
\begin{equation}\label{aml}
\frac{a_m(l)}{c_m(l)} \approx b_m(l),
\end{equation}
since~\(0 \leq a_m(l) \leq 1+\frac{D}{m}\) and~\(c_m(l) \geq 1-\frac{D}{m}\) for all~\(m\) large and some constant~\(D > 0\) and so~\[\left|\frac{a_m(l)}{c_m(l)} - a_m(l)\right| = \frac{a_m(l)}{c_m(l)} \left|1-c_m(l)\right| \leq \left(\frac{1+Dm^{-1}}{1-Dm^{-1}}\right)\frac{D}{m} \leq \frac{2D}{m}\]
for all~\(m\) large. Thus~\(\frac{a_m(l)}{c_m(l)} \approx a_m(l) \approx b_m(l).\) We now write
\begin{equation}\label{gm_l_eval}
G_m(l) = \sum_{l_1=0}^{l-1} g_m(l_1) = \frac{\frac{1}{m}\sum_{l_1=0}^{l-1}f\left(\frac{l_1}{m}\right)}{\frac{1}{m}\sum_{l_1=0}^{Mm}f\left(\frac{l_1}{m}\right)} =: \frac{a_m(l)}{c_m(l)}.
\end{equation}
Using~(\ref{gapp2}), we have that
\[a_m(l) = \frac{1}{m}\sum_{l_1=0}^{l-1}f\left(\frac{l_1}{m}\right) \approx F\left(\frac{l}{m}\right) =: b_m(l)\] and  using the first relation in~(\ref{gapp4}), we get~\(c_m(l) = \frac{1}{m}\sum_{l_1=0}^{Mm}f\left(\frac{l_1}{m}\right) \approx 1.\) From~(\ref{aml}), we then get the second relation in~(\ref{gapp4}).

To prove~(\ref{gapp3}) write
\[A_m(l) = \frac{1}{m}\sum_{l_1=0}^{l-1} F^{r-1}\left(\frac{l_1}{m}\right) f\left(\frac{l_1}{m}\right) \leq \frac{1}{m} \sum_{l_1=0}^{Mm} f\left(\frac{l_1}{m}\right)  \approx 1\]
and~\(C_m(l) = \frac{1}{m}\sum_{l_1=0}^{Mm}f\left(\frac{l_1}{m}\right) \approx 1,\) by the first relation in~(\ref{gapp4}). Arguing as in~(\ref{aml}), we have that~\(\frac{A_m(l)}{C_m(l)} \approx A_m(l)\) and from~(\ref{gapp2}), we have that~\(A_m(l) \approx \frac{1}{r}F^{r}\left(\frac{l}{m}\right),\) proving~(\ref{gapp3}).~\(\qed\)

\emph{Proof of Theorem~\ref{t1t2_thm}}: For simplicity, we let~\(X_{j} = X_{m,j}\) throughout
and define~\(B_k := A_{i_1} \bigcap \ldots \bigcap A_{i_k}.\)
As before, we prove a slightly stronger result. For~\(l > 0,\) we let
\[A_{i_k}(l) = A_{i_k} \bigcap \left\{X_{i_k} < \frac{l}{m} \right\}\] and
\begin{equation}\label{bkl_def}
B_k(l) := \bigcap_{j=1}^{k-1} A_{i_j} \bigcap A_{i_k}(l)
\end{equation}
and prove that
\begin{equation}\label{asym2_disc}
\mathbb{P}_m(B_k(l)) \approx \frac{F^{i_k}\left(\frac{l}{m}\right)}{c(i_1)\cdots c(i_k)}.
\end{equation}

In the proof of~(\ref{asym2_disc}) below, we use the following fact throughout:\\
If~\(p_m(l) \approx q_m(l),\) then using~\(\sum_{l_1=0}^{l-1} g_m(l_1) \leq 1\) we also have that
\begin{equation}\label{am_rel_eq}
\sum_{l_1=0}^{l-1} p_m(l_1) g_m(l_1) \approx \sum_{l_1=0}^{l-1} q_m(l_1) g_m(l_1).
\end{equation}

\emph{Proof of~(\ref{asym2_disc}) for~\(k=1\)}: Letting~\(n \geq 1\) be fixed we have for~\(1 \leq k \leq n \) that
\begin{eqnarray}
\mathbb{P}_m\left(A_{i_k}(l)\right) &=& \mathbb{P}_m\left(\max_{j \in {\cal C}(i_k)} X_j < X_{i_k} < \frac{l}{m}\right) \nonumber\\
&=& \sum_{l_1=0}^{l-1} \mathbb{P}_m\left(\max_{j \in {\cal C}(i_k)}  X_j < \frac{l_1}{m}\right) g_m(l_1)\nonumber\\
&=& \sum_{l_1=0}^{l-1} G_m^{c(i_k)-1}(l_1) g_m(l_1) \nonumber\\
&\approx& \sum_{l_1=0}^{l-1} F^{c(i_k)-1}\left(\frac{l_1}{m}\right) g_m(l_1) \label{x1}\\
&\approx& \frac{1}{c(i_k)} F^{c(i_k)}\left(\frac{l}{m}\right),\label{pr_ak_ev1_u}
\end{eqnarray}
where~(\ref{x1}) follows from~(\ref{am_rel_eq}) and the fact that~\(G_m^{r}(l_1) \approx F^{r}\left(\frac{l_1}{m}\right)\) (see~(\ref{gapp2})) and~(\ref{pr_ak_ev1_u}) follows using~(\ref{gapp3}) in Lemma~\ref{thet_lem}.

To prove the induction step, we now assume that the relation~(\ref{asym2_disc}) is true some integer~\((k-1) \geq 2.\) If the event~\(B_{k-1} := \bigcap_{j=1}^{k-1} A_{i_j}\) occurs, then arguing as in~(\ref{ai2_ai1_k}), we get that~\(\ind(B_k(l)) = \ind\left(A_{i_k}(l) \bigcap B_{k-1}\right)\) equals
\begin{equation}
\ind\left(\frac{l}{m} > X_{i_k} > \max_{j \in {\cal C}(i_k) \setminus \left({\cal C}(i_{k-1}) \cup i_{k-1}\right)} X_j\right) \cdot \ind\left(\frac{l}{m} > X_{i_k} > X_{i_{k-1}}\right) \ind(B_{k-1}). \label{ai2_ai1_k_disc}
\end{equation}
Again, the event~\( \ind(B_{k-1})\) depends only on the values of~\(\{X_j\}_{1 \leq j \leq i_{k-1}}\) and so conditioning on~\(X_{i_k}= \frac{l_1}{m}\)
we get from~(\ref{ai2_ai1_k}) that
\begin{equation}\label{pa2a1_disc}
\mathbb{P}\left(B_k(l)\right) = \sum_{l_1 = 0}^{l-1} t_{i_1,i_2,\ldots,i_k}(l_1) g_m(l_1),
\end{equation}
where~
\begin{eqnarray}
t_{i_1,i_2,\ldots,i_k}(l_1) &=&  \mathbb{P}_m\left(\max_{j \in {\cal C}(i_k) \setminus \left({\cal C}(i_{k-1}) \cup i_{k-1}\right)} X_j < \frac{l_1}{m}\right) \mathbb{P}\left(B_{k-1} \bigcap \left\{ X_{i_{k-1}} < \frac{l_1}{m}\right\} \right) \nonumber\\
&=& G_m^{c(i_k)-c(i_{k-1})-1}(l_1)\mathbb{P}\left(B_{k-1} \bigcap \left\{ X_{i_{k-1}} < \frac{l_1}{m}\right\} \right) \nonumber\\
&=& G_m^{c(i_k)-c(i_{k-1})-1}(l_1)\mathbb{P}\left(B_{k-1}(l_1)\right) \nonumber
\end{eqnarray}
so that
\begin{eqnarray}
\mathbb{P}\left(B_k(l)\right) &=& \sum_{l_1=0}^{l-1} G_m^{c(i_k)-c(i_{k-1})-1}(l_1)\mathbb{P}\left(B_{k-1}(l_1)\right) g_m(l_1) \nonumber\\
&\approx& \sum_{l_1=0}^{l-1} F^{c(i_k)-c(i_{k-1})-1}\left(\frac{l_1}{m}\right)\mathbb{P}\left(B_{k-1}(l_1)\right) g_m(l_1) \label{ty1}\\
&\approx& \sum_{l_1=0}^{l-1} F^{c(i_k)-c(i_{k-1})-1}\left(\frac{l_1}{m}\right)\left(\frac{F^{c(i_{k-1})}\left(\frac{l_1}{m}\right)}{c(i_1)\cdots c(i_{k-1})}\right) g_m(l_1).\;\; \label{ty2}
\end{eqnarray}
The expression~(\ref{ty1}) follows from the fact that~\(G^{r}_m(l) \approx F^{r}\left(\frac{l}{m}\right)\) (see~(\ref{gapp2})) and the summation relation~(\ref{am_rel_eq})
and~(\ref{ty2}) follows again from the summation relation~(\ref{am_rel_eq}) and the induction assumption. We recall from the strictly increasing comparison set sequence condition~\((a1)\) that~\(c(i_k) \geq c(i_{k-1})+1.\)

From~(\ref{ty2}) we get that
\begin{eqnarray}
\mathbb{P}\left(B_k(l)\right) &\approx& \left(\frac{1}{c(i_1)\cdots c(i_{k-1})}\right) \sum_{l_1=0}^{l-1} F^{c(i_k)-1}\left(\frac{l_1}{m}\right)g_m(l_1) \nonumber\\
&\approx&  \frac{1}{c(i_1)\cdots c(i_{k})} \cdot F^{c(i_k)}\left(\frac{l}{m}\right) \nonumber
\end{eqnarray}
using the estimate~(\ref{gapp3}) from Lemma~\ref{thet_lem}, proving the induction step.~\(\qed\)




\subsection*{Acknowledgement}
I thank Professors Alberto Gandolfi and Federico Camia for crucial comments and for my fellowships.

\bibliographystyle{plain}

\end{document}